# ON THE TWO-PHASE FRAMEWORK FOR JOINT MODEL AND DESIGN-BASED INFERENCE

By Susana Rubin-Bleuer[1] and Ioana Schiopu Kratina

*Statistics Canada*

We establish a mathematical framework that formally validates the two-phase "super-population viewpoint" proposed by Hartley and Sielken [*Biometrics* **31** (1975) 411–422] by defining a product probability space which includes both the design space and the model space. The methodology we develop combines finite population sampling theory and the classical theory of infinite population sampling to account for the underlying processes that produce the data under a unified approach. Our key results are the following: first, if the sample estimators converge in the design law and the model statistics converge in the model, then, under certain conditions, they are asymptotically independent, and they converge jointly in the product space; second, the sample estimating equation estimator is asymptotically normal around a super-population parameter.

**1. Introduction.** Classical sampling theory concerns inference for finite population parameters. For the finite population mean $\overline{Y} = \sum_{i=1}^{N} y_i/N$, inference typically considers the interval

$$[\bar{y} \pm t_p se(\bar{y})],$$

where $t_p$ is a constant chosen with a normal or Student distribution in mind, and $se(\bar{y})$ denotes the standard error of the sample mean $\bar{y}$ (see [13]). The expression above means that $\overline{Y}$ is within the interval $[\bar{y} - t_p se(\bar{y}), \bar{y} + t_p se(\bar{y})]$ with some degree of confidence. Here $N$ is the size of the finite population, the $y_i$'s are considered nonstochastic but unknown numbers and probability statements arise from the selection of units in the sample. No distributional assumptions are made about the $y_i$'s. This nonparametric approach to inference is often called design-based inference.

Received January 2003; revised December 2004.
[1]Supported in part by the Sabbatical Research Program of Statistics Canada.

However, there are many situations when we have to resort to postulating a model. For descriptive analysis in a finite population, we need a model when we have to deal with nonresponse, small area estimation or measurement errors. For studies involving scientific questions, the parameters of stochastic models are sometimes of more interest than finite population parameters. For example, in longitudinal surveys we are interested in modeling the dependencies between health status and certain socio-economic covariates. Staying within the finite population framework limits our inference to the reference population only. To illustrate the issue, let us take the sample mean $\bar{y}$ obtained from a sample of size $n$, and suppose that we wish to draw conclusions on a more general population than the finite population from which we obtained the sample: we view $\bar{y}$ as an estimator of the model mean $\mu$. We have

$$(1.1) \qquad \sqrt{n}(\bar{y} - \mu) = \sqrt{n}(\bar{y} - \overline{Y}) + \sqrt{(n/N)}\sqrt{N}(\overline{Y} - \mu).$$

The design-based large sample properties of the first term on the right-hand side of (1.1) have been studied for many sampling designs. Conditions were given for the asymptotic normality of the sample mean (design-based central limit theorem, or CLT): for simple random sampling without replacement (SRSWOR) and rejective sampling with varying probabilities by Hájek [8, 9], for probability proportional to size without replacement ($\pi ps$) designs by Rosén [19, 20], and for stratified multistage probability proportional to size with replacement (PPSWR) designs by Krewski and Rao [15]. For descriptions of these and other sampling designs, see, for example, [28]. To derive a design-based CLT for the left-hand side of (1.1), we would have to assume not only that the sampling rate $n/N$ converges to zero, but also that the sequence of numbers $\sqrt{N}(\overline{Y} - \mu)$ is bounded as $N \to \infty$. As a sequence of numbers, this last condition is very restrictive. However, as a sequence of sums of independent, identically distributed (i.i.d.) random variables (r.v.) in the super-population, $\sqrt{N}(\overline{Y} - \mu)$ is bounded in probability and the second term of the right-hand side of (1.1) converges to zero in the probability of the model when $n/N$ converges to zero. To study the asymptotic properties of the survey sample means around the model mean, it is necessary to include the model and the design in the same probability space.

In this article we first construct a product space, which is a mathematical framework for joint design-based and model-based inference. Our key results are quite general. First, we show that, under certain conditions, if the survey sample estimators converge in the law of the sampling design and the associated model statistics converge in the law of the model (not necessarily to a Gaussian distribution), then they are asymptotically independent and they converge jointly in the product space. Second, we show that a survey sample estimator of a model parameter, which is derived from a very general sampling estimating equation, exists, is consistent and is asymptotically normal. Hartley and Sielken [11] introduced the "super-population" approach to describe the relationship between the infinite population (also called superpopulation) and the finite population from which we select the sample. Many authors worked within the two-phase framework and accounted for the variability due to the design and the model by means of the "anticipated variance." The contributions of Fuller [6], Isaki and Fuller [12], Godambe and Thompson [7], Korn and Graubard [14], Pfeffermann and Sverchkov [17], Binder and Roberts [4], Rodríguez [18] and Molina, Smith and Sugden [16] are just a few among the many on the subject.

Fuller [6] established large sample properties of the sample regression estimator around the model parameter with data obtained from stratified cluster samples. His approach could only be applied to stratified SRSWOR designs in the first stage. Our general approach to estimation of model parameters extends Fuller [6] to more general designs and estimators, even if the sampling rate is nonnegligible.

The formal expression of the product space, together with the key results described above, establish a general and unified methodology that accommodates the diverse techniques of these authors, and has enabled us to extend some of their results. Moreover, the more formal aspects of the methodology (sub-$\sigma$-fields, filtrations in the product space) proved to be essential for adapting counting process methodology to the analysis of survival survey data (see [23, 24]). In addition, the design-based distribution of a sample estimator is a "second phase" concept, that is, a conditional distribution given the minimal information in the model. In general, we could apply this methodology to most situations where we have a two phase randomization process.

The joint design-model distribution of the sample data is also called the distribution of the sample variables (see [17]). We present other results that refer to the sample variables under the posterior distribution given the sample labels. For a sequence of random variables $Y = \{Y_1, \ldots, Y_N\}$, it is well known that the posterior distribution of $Y$ given the sample outcome $\{(i, Y_i = y_i), i \in s_0\}$ depends only on the sample $s_0$ actually drawn and not on the sampling design used to draw it, provided that $s_0$ and $Y$ are stochastically independent given the design variables (see [29]). Here, however, we look at conditioning just on the $s_0$ actually drawn. We also show that the posterior distribution depends only on $s_0$. Note that whether the labels are repeated or not is a consequence of the design. If the sample $s_0$ from a with replacement (WR) sampling design has repeated labels, the sample variables under the posterior distribution given $s_0$ are not stochastically independent even if the original components of $Y$ were independent. It is also well known that an SRSWOR from a finite population, which was generated by a super-population, when viewed as a sample from the infinite space inherits the same properties of the random variables which generated the finite population. Fuller [6] applied the CLT to the array of variables from an SRSWOR design to obtain the asymptotic distribution of the sample regression estimator. It is not clear why the classical CLT could be applied in [6] without further assumptions. In this article, we show formally that his array of sample variables, not necessarily nested, consists of i.i.d. variables under the posterior distribution given the sample labels. For the CLT to hold for this array, we only require that the original super-population variables have a finite variance.

In order to obtain the total (anticipated) variance in (1.1), we must impose (model-based) conditions on the super-population model, which survey statisticians would rather avoid. At the very least, some form of model-based independence is needed. Many authors assume that the sampling rate is small enough so they can ignore the variation due to the model component. However, Korn and Graubard [14] show that we should not dismiss the second term in the total variance without checking first that it is indeed sufficiently small relative to the first term.

The article is organized as follows. In Sections 2–5 we develop the tools needed to incorporate the design and the model in the same space. Section 6 is an application of the product space methodology. In Section 2 we modify somewhat the usual definitions of sample design, finite population parameter and sample estimator to enable us to view them as random variables in the super-population (Definition 4.2, Remark 4.2). In Section 3 we adopt the super-population definition in [28] to define what it means for a finite population to be generated by a super-population (Definition 3.1). Proposition 3.1 shows how conditions needed for the design-based CLT follow from simple conditions in the super-population. In Section 4 we define the general product space (Definitions 4.1, 4.3) and show how stochastic dependence is introduced in the product space (Example 4.1). We exploit the additional information on the design and the model by calculating posterior distributions and we study the interplay between dependence and independence of random variables viewed in the design space, the product space or the model space (Example 4.1 and Proposition 4.2). In Section 5 we show that, if the sample and super-population statistics converge in law in their respective spaces, they also converge in law in the product space. The two terms in the right-hand side of (1.1) are not, in general, stochastically independent. We establish here their "asymptotic independence" under mild conditions in Theorem 5.1. Example 5.1 yields the asymptotic normality of the ratio estimator of the weighted average of the strata means under a stratified one-stage PPSWR design. In Section 6 we establish the existence and asymptotic normality of a sample estimator derived from a general estimating equation, under general conditions. Example 6.1 is an application to a two-stage sampling design.

## 2. Finite populations and sampling designs.

DEFINITION 2.1. A finite population $U = \{1, \ldots, N\}$ of size $N$ consists of $N$ labels, with their associated data, that is, each unit $i$ is associated to a unique vector $(y_i, x_i, z_i)$, $i = 1, \ldots, N$. Here $y_i \in \mathbb{R}^p$, $x_i \in \mathbb{R}^k$ represent, respectively, the characteristics of interest and the auxiliary information, and $z_i \in \mathbb{R}_+^q$ is the "prior" information available at the time of the design of the survey on all units $i = 1, \ldots, N$. We write $y^N = (y_i)_{i=1,\ldots,N}$, $x^N = (x_i)_{i=1,\ldots,N}$ and $z^N = (z_i)_{i=1,\ldots,N}$.

REMARK 2.1. In this paper $N$ will denote the size of the finite population (i.e., the number of ultimate sampling units in the population) for one-stage-sampling schemes, and it will denote the number of clusters or primary sampling units (p.s.u.s) for multistage schemes, in which case the size of the finite population will be denoted by $M$.

DEFINITION 2.2. A sample is the realization of a probabilistic (randomized) selection or sampling scheme ([28], page 25). We adopt the comprehensive definition of a sample in [10], page 42: it views the sample as "a finite sequence of units or labels of the finite population, which are drawn one by one until the sampling is finished according to some stopping rule. This sequence distinguishes the order of units, may be of variable length and may include one unit of the finite population several times." This definition includes samples selected without replacement (WOR) and WR. In what follows, we do not require that samples be selected sequentially, but, for convenience, we may consider an order in which the $n$ sampled units are either observed or selected.

In the literature, a design $p$ associated with a sampling scheme is a probability function on the set of all possible samples under this scheme (see, e.g., [28]). The definition of a sampling design given below requires measurability of $p$ as a function of the variables containing the prior information. The same holds for Definition 2.4 of a finite population parameter (cf. [28], page 39). The measurability conditions ensure that, when the finite population is generated by a super-population, the finite population parameter and the estimator are real-valued measurable functions (random variables) defined on the probability space associated with the super-population (see Definition 3.1).

DEFINITION 2.3. Let $U$ be the finite population of Definition 2.1. Given a sampling scheme, let $S$ be the set of all possible samples under the scheme. Let $C(S)$ consist of all subsets of $S$. A sampling design associated to a sampling scheme is a function $p : C(S) \times \mathbb{R}_+^{q \times N} \to [0, 1]$ such that:

(i) for all $s$ in $S$, $p(s, \cdot)$ is Borel-measurable in $\mathbb{R}_+^{q \times N}$;
(ii) for $z^N \in \mathbb{R}_+^{q \times N}$, $p(\cdot, z^N)$ is a probability measure on $C(S)$.

We say that $(S, C(S), p)$ is a design probability space, where $p(s, \cdot) > 0$, $s \in S$.

REMARK 2.2. For the sake of simplicity, in all applications we will take $q = 1$.

REMARK 2.3. For a one-stage Poisson sampling scheme, the collection $S$ of all possible samples is completely determined given only the sizes of the strata in the population. For other sampling schemes, $S$ cannot be determined unless we know the strata and sample sizes, and possibly other parameters, depending on the sampling scheme. Under a one-stage $\pi ps$ scheme, $S$ can be defined without prior knowledge of the unit sizes. Under a first-stage $\pi ps$ scheme and a second stage SRSWOR scheme, we cannot completely determine $S$ unless we know the first-stage unit sizes, since they are the second-stage population sizes.

DEFINITION 2.4. Consider a finite population as in Definition 2.1. A finite population parameter $\theta_N$ is a Borel-measurable function defined on a subset of $\mathbb{R}^{(p+k+q) \times N}$. An estimator of this finite population parameter associated with a design, also called a sample estimator, is a function $\hat{\theta}_N : S \times \mathbb{R}^{(p+k+q) \times N} \to \mathbb{R}$, where $\hat{\theta}_N(s, \cdot)$ is Borel-measurable.

In the next example we define the fundamental notation used by Krewski and Rao [15], which we will use subsequently (e.g., in Proposition 3.1).

EXAMPLE 2.1 (Stratified two-stage PPSWR [15]). Let $N_h$ be the number of p.s.u.s in stratum $h$, $M_{hi}$ be the number of ultimate units in p.s.u. $hi$, $i = 1, \ldots, N_h$, $h = 1, \ldots, L$, and $L$ the number of strata. Let $N = \sum_{h=1}^{L} N_h$, $M_h = \sum_{i=1}^{N_h} M_{hi}$ and $M = \sum_{h=1}^{L} M_h$. The prior information consists of the "sizes" $z_{hi} = M_{hi}$, $i = 1, \ldots, N_h$, $h = 1, \ldots, L$. Suppose $n_h \geq 2$ p.s.u.s are selected with replacement in stratum $h$ with probabilities $p_{hi} = M_{hi}/M_h$, $i = 1, \ldots, N_h$, $h = 1, \ldots, L$ at each draw. The selection is independent in each stratum, and independent second stage samples are taken within those p.s.u.s selected more than once. The finite population mean is $\theta_N = \sum_{h=1}^{L} W_h \theta_h$, where $W_h = M_h/M$ is the stratum weight, $\theta_h = \overline{Y}_h = \sum_{i=1}^{N_h} y_{hi}/M_h$ is the finite population stratum mean and $y_{hi}$ is the total of p.s.u. $hi$, $i = 1, \ldots, N_h$, $h = 1, \ldots, L$. Let $I_{hi}^k = 1$ if p.s.u. $hi$ is selected in the sample at the $k$th draw in stratum $h$ and 0 otherwise, $k = 1, \ldots, n_h$, $i = 1, \ldots, N_h$, $h = 1, \ldots, L$. If the cluster $hi$ is selected at the $k$th draw, let $\hat{y}_{hi}$ be an unbiased estimator of the total $y_{hi}$ based on sampling at the second stage and set $\hat{y}_{hi} = 0$ otherwise, $k = 1, \ldots, n_h$, $i = 1, \ldots, N_h$, $h = 1, \ldots, L$. For stratum $h$, we consider the estimator $\hat{\theta}_h = \sum_{k=1}^{n_h} \hat{\theta}_h^k / n_h$, where $\hat{\theta}_h^k = \sum_{i=1}^{N_h} \hat{y}_{hi} I_{hi}^k / M_{hi}$, $k = 1, \ldots, n_h$, $i = 1, \ldots, N_h$, $h = 1, \ldots, L$. Finally, a design-unbiased sample estimator of $\theta_N$ is $\hat{\theta}_N(y^N, M^N) = \sum_{h=1}^{L} W_h \hat{\theta}_h$.

We often refer to conditions $C_1$ to $C_3$ of Yung and Rao [31], which evolved from conditions introduced by Krewski and Rao [15] for the asymptotic normality of the sample mean $\hat{\theta}_N$ (see the Appendix).

## 3. Super-populations.

DEFINITION 3.1. Consider a finite population $U$ of size $N$ as in Definition 2.1. A super-population associated with it consists of a probability space $(\Omega, \mathscr{F}, P)$ and random vectors $(Y_i, X_i, Z_i)$, $Y_i : \Omega \to \mathbb{R}^p$, $X_i : \Omega \to \mathbb{R}^k$, $Z_i : \Omega \to \mathbb{R}^q_+$, such that $Y_i(\omega_0) = y_i$, $X_i(\omega_0) = x_i$, $Z_i(\omega_0) = z_i$, for some $\omega_0 \in \Omega$, $i = 1, \ldots, N$. We write $Y^N = (y_i)_{i=1,\ldots,N}$ and define $X^N$ and $Z^N$ similarly. We say that $U$ is a realization of or is generated by the super-population. Any distribution of $(Y^N, X^N, Z^N)$ that is given a priori is called a super-population model. We note that different outcomes $\omega$ can generate the same finite population.

Definition 3.1 is similar to the definition given in [28], page 533. We assume throughout this work that $N$ is not random. In what follows, the subscript "$d$" refers to design randomization and "$m$" refers to the randomization on the probability space $(\Omega, \mathscr{F}, P)$. We use $E_m$, $V_m$ to denote, respectively, the expectation and variance with respect to the probability space $(\Omega, \mathscr{F}, P)$. We use the standard notation $\sigma(X)$ for the $\sigma$-field generated by the function $X$ (see also Definition 3.1 in [22] or [25]).

EXAMPLE 3.1 (Two-stage super-population model). Let $\Omega$ be the conceptual population of people living in a country. Suppose it is composed of $L$ disjoint strata of units $hi$, $i = 1, \ldots, N_h$, $h = 1, \ldots, L$, where unit $hi$ represents a cluster of individuals. Let $(\Omega, \mathscr{F}, P)$ be the corresponding probability space. Now we assume that $Z_{hi}$ are discrete r.v.s on the probability space that represent the number of individuals that live in cluster $hi$. We are interested in characteristics $Y_{hij}$ pertaining to the individuals labelled by $hij$, living in cluster $hi$, $i = 1, \ldots, N_h$, $h = 1, \ldots, L$. In order to be able to define the super-population according to Definition 3.1, we must know an outcome of the $Z_{hi}$, say, the sizes of the clusters of the population existing right now. Let $F_M = \{\omega \in \Omega : Z_{hi}(\omega) = M_{hi}, i = 1, \ldots, N_h, h = 1, \ldots, L\}$. We use this information to define the super-population model by conditioning on the $\sigma$-field generated by the event $F_M$. The conditional probability measure is given by $P_M(F, \omega_0) = P(F|F_M)$ if $\omega_0 \in F_M$ for $F \in \mathscr{F}$ (see [5], Equation 3, page 222 and note that $P(F_M) > 0$ since the r.v.s $Z_{hi}$ are discrete). Now we define the super-population on $(\Omega, \mathscr{F}, P_M)$ by random vectors $Y_{hij}$ of $p$ socio-economic characteristics associated with the individual $hij$, $Y_{hij} : \Omega \to \mathbb{R}^p$, $j = 1, \ldots, M_{hi}$, $i = 1, \ldots, N_h$, $h = 1, \ldots, L$. The cluster totals $\{Y_{hi} = \sum_{j=1}^{M_{hi}} Y_{hij}, i = 1, \ldots, N_h, h = 1, \ldots, L\}$ are assumed i.i.d. r.v.s within strata.

We now illustrate how conditions that are sufficient for design-based CLTs can be justified as a consequence of simple moment conditions in the super-population, which, in turn, can be justified by expert knowledge of the model.

Consider the two-stage super-population model of Example 3.1 and assume that the total number of clusters $N \to \infty$. Assume the sampling design of Example 2.1, defined on the finite population generated by $\omega \in \Omega$, where $Y_{hi}(\omega) = \sum_{j=1}^{M_{hi}} Y_{hij}(\omega)$, $\omega \in \Omega$. In Proposition 3.1 below we show that moment conditions in the super-population yield the Liapunov-type condition $(C_1')$ similar to $\sum_{h=1}^{L} W_h E_d |\hat{\theta}_h^k - \overline{Y}_h|^{2+\delta} = O(1)$ as $n \to \infty$, $\hat{\theta}_h^k$ as in Example 2.1, which is condition $C_1$ of Krewski and Rao [15].

PROPOSITION 3.1. *Let $n = n_1 + n_2 + \cdots + n_L$. We assume the model-based condition*

$$(M_1) \qquad (1/N) \sum_{h=1}^{L} \sum_{i=1}^{N_h} E_m |Y_{hi}|^{2+\delta} = O(1), \qquad \delta > 0 \text{ as } N \to \infty.$$

*Then*

$(C_1')$ *for all $k = 1, \ldots, n_h$, $h = 1, \ldots, L$, $\sum_{h=1}^{L} W_h E_d |\hat{\theta}_h^k - \overline{Y}_h|^{2+\delta}(\omega) = O(1)$,*

*for all $\omega$ in a set with model probability 1 (a.s. $\omega$), $N \to \infty$, where $\hat{\theta}_h^k$ is the estimator of the stratum mean based on the $k$th draw in stratum $h$, $1 \leq k \leq n_h$, $h = 1, \ldots, L$, defined in Example 2.1.*

The proof is given in the Appendix. Here $E_d$ is the design-based expectation and is calculated in the Appendix. Note that $E_d |\hat{\theta}_h^k|^{2+\delta}(\omega)$ is a random variable in the model space.

**4. The product space.** In this section we define a product probability space that includes the super-population and the design space, under the premise that sample selection and the model characteristic $Y$ are independent given all of the design variables $Z$. We investigate independence properties of the sample variables under the posterior distribution given the sample labels and we provide the formal proof of the CLT under the posterior distribution for an SRSWOR design. Proposition 4.4 derives the product space probability given the model.

DEFINITION 4.1. Consider a finite population of size $N$ generated by a super-population $(Y^N, X^N, Z^N)$ as in Definition 3.1. We define the product space as the set $S \times \Omega$ with the $\sigma$-field $C(S) \times \mathscr{F}$.

DEFINITION 4.2. Consider a super-population associated with a finite population as in Definition 4.1. Let $p: C(S) \times \mathbb{R}_+^{q \times N} \to [0, 1]$ be a sampling design on the finite population as in Definition 2.3. Then the sampling design can be viewed as a random variable on $(S \times \Omega, C(S) \times \mathscr{F})$ defined by

$$(4.1) \qquad p(s, \omega) = p(s, Z_N(\omega)), \qquad s \in S, \ \omega \in \Omega.$$

DEFINITION 4.3. We define $P_{d,m}$ as the $\sigma$-additive measure that, on elementary rectangles of the product $\sigma$-field, has the value

$$(4.2) \qquad P_{d,m}(\{s\} \times F) = \int_F p(s, \omega)\, dP, \qquad s \in S, \ F \in \mathscr{F}.$$

Note that each set in $C(S) \times \mathscr{F}$ can be expressed as a finite union of elementary rectangles and $P_{d,m}(S \times \Omega) = 1$. Hence, $P_{d,m}$ is a probability measure on the product space. If $Z^N$ are discrete random variables, we may build the product space from the super-population model given $Z^N$, with the probability measure $P_z(\cdot) = P(\cdot|F_z)$, $F_z = \{\omega : Z^N(\omega) = z^N\}$. With $P_z$ replacing $P$ in (4.2), we obtain

$$P_{d,m}^*(\{s\} \times F) = p(s, z^N) \cdot P_z(F), \qquad s \in S, \ F \in \mathscr{F}.$$

REMARK 4.1. Any measurable set in the product space is of the form $B = \bigcup_{s \in S} \{s\} \times F_s$, where some sets $F_s \in \mathscr{F}$ could be empty. By Definition 4.3, $P_{d,m}(B) = \int_\Omega \sum_{s \in S} p(s, \omega) I_{F_s}(\omega)\, dP$. We denote the integrand by $P_{d,m}(B|S \times \mathscr{F})(\omega) = \sum_{s \in S} p(s, \omega) I_{F_s}(\omega)$. By Proposition 4.4 this is a conditional probability given the $\sigma$-field $S \times \mathscr{F}$.

REMARK 4.2. Let $\hat{\theta}_N$ be a sample estimator on the design space with associated super-population $(Y^N, X^N, Z^N)$. It can be viewed as a random variable on the product space defined by

$$(4.3) \quad \hat{\theta}_N(s, \omega) = \hat{\theta}_N(s, Y^N(\omega), X^N(\omega), Z^N(\omega)), \qquad s \in S, \ \omega \in \Omega.$$

We omit writing the index $N$ when no confusion may arise.

DEFINITION 4.4 (The sample variables). The components of a sample outcome $y^s = \{y_i, i \in s\}$, $s \in S$, can be viewed as random variables in the product space, and following Pfeffermann and Sverchkov [17], we call them sample variables.

A sample can be written as a sequence of labels $i(k)$, indexed by $k = 1, \ldots, n$, the order in which the labels are observed. Let us define $I_i^k = 1$ if label $i(k) = i$ and $I_i^k = 0$ otherwise. If the sample is drawn sequentially, the $I_i^k$ coincide with the $k$th draw indicators in Example 2.1. Thus, the sample outcome can be written as the sequence of $n$ units, where each coordinate $k$ of the sequence represents the $y$-value for the label $i(k)$, $k = 1, \ldots, n$:

$$y^s = \left( \sum_{i=1}^N y_i I_i^1(s), \sum_{i=1}^N y_i I_i^2(s), \ldots, \sum_{i=1}^N y_i I_i^n(s) \right).$$

The *sample variables* can be written as

$$Y_{i(k)}(s, \omega) = \sum_{j=1}^N Y_j(\omega) I_j^k(s), \qquad k = 1, \ldots, n.$$

We will use this notation subsequently. Note that, for WR designs, the labels $i(k)$ and $i(l)$ could be the same for $k \neq l$.

REMARK 4.3. Assume that the components of $Y^N$ are independent random variables. If the design is SRSWOR and the components of $Y^N$ are i.i.d. in the super-population, the "sample variables" $Y_{i(k)}$, $k = 1, \ldots, n$, are independent in the product space. However, if the original $Y^N$ are not identically distributed, the variables $Y_{i(k)}$, $k = 1, \ldots, n$, may become stochastically dependent in the product space. Under a simple random sample with replacement (SRSWR) design, the variables $Y_{i(k)}$, $k = 1, \ldots, n$, are stochastically dependent in the product space whether the original super-population variables are i.i.d. or not. We refer to the Appendix for an illustration of the mechanism.

EXAMPLE 4.1 (Stochastic dependence in the product space). Let $N = n = 2$ and $Y_i$, $i = 1, 2$, be i.i.d. r.v.s each with a Bernoulli distribution $B(1, 0.5)$. Under simple random sampling (SRS), $P_{d,m}(Y_{i(1)} = 1) = P_{d,m}(Y_{i(2)} = 1) = 0.5$. Under SRSWR, $P_{d,m}(Y_{i(1)} = 1, Y_{i(2)} = 0) = 0.125 \neq 0.5 \times 0.5$ [see (A.2′) in the Appendix], whereas under SRSWOR, $P_{d,m}(Y_{i(1)} = 1, Y_{i(2)} = 0) = 0.25 = 0.5 \times 0.5$.

EXAMPLE 4.2 (Two-stage super-population model and two stage design). We assume the two-stage super-population model of Example 3.1, where we use the size of the clusters of a population existing right now to define the model. This minimum necessary information is contained in $F_M = \{ \omega \in \Omega : Z_{hi}(\omega) = M_{hi}, i = 1, \ldots, N_h, h = 1, \ldots, L \}$, where the $M_{hi}$ are cluster sizes as in Example 3.1. We select the sample with probability proportional to those sizes, but we want to draw conclusions about a more general population than the finite population living in those clusters now. We set $M^{N_h} = (M_{hi})_{i=1,\ldots,N_h}$, $h = 1,\ldots,L$. Once the model is defined, we define a sample space $S$ as the collection of all possible "stratified clustered" sequences of units (see Remark 2.3) of a finite population associated with the super-population model. Then we define a stratified two-stage sampling design $p(s, M^{N_1}, \ldots, M^{N_L})$ with $L$ strata, $N$ clusters and $M$ ultimate units. We then construct the space $S \times \Omega$ with probability measure $P_{d,m}$ defined on the elementary rectangles by $P_{d,m}(s \times F) = p(s, M^{N_1}, \ldots, M^{N_L})P_M(F)$, $s \in S$, $F \in \mathscr{F}$ (see also Example 4.3 in [27]).

Consider a sample $s_0 \in S$ and let $\sigma(s_0 \times \Omega)$ be the four-set sub-field generated by $s_0 \times \Omega$. Let $P(\cdot|s_0)$ be the conditional probability measure given this field. We have the following result.

PROPOSITION 4.1. *For each $B = \bigcup_{s \in A}\{s\} \times F_s$, $A \in C(S)$, $F_s \in \mathscr{F}$, we have*

(4.4)      (i)     $P(B|s_0) = P_{d,m}(s_0 \times F_{s_0})/P_{d,m}(s_0 \times \Omega)$

*if $s_0 \in A$, and $0$ otherwise. If, in particular, $p(s, \omega)$ does not depend on $\omega \in \Omega$, and $I_A(s_0)$ is the value of the indicator function of the set $A$ at $s_0$, we have*

(ii)     $P(B|s_0) = P(F_{s_0})I_A(s_0).$

PROOF. (i) is immediate from [5], Example 1, page 223. Statement (ii) follows from (i). □

PROPOSITION 4.2 [Stochastic independence of the sample under $P(\cdot|s_0)$]. *Let $Y^N$ denote the super-population composed of $N$ independent random vectors. Assume an SRS design. Under $P(\cdot|s_0)$, the $Y_{i(k)}$, $k = 1,\ldots,n$, variables are stochastically independent if there are no repeated labels in the selected sample and stochastically dependent otherwise.*

See the Appendix for the proof (see also [26]).

EXAMPLE 4.3. Let $N = n = 2$. Suppose, as in Example 4.1, that $Y_i$, $i = 1,2$, are i.i.d. r.v.s distributed as $B(1, 0.5)$. Assume that we selected $s_0 = \{1,2\}$ under SRS. This sample has no repeated labels and $I_i^1(s_0)I_j^2(s_0) = 0$ if $i = j$, $i,j = 1,2$. Then $P(Y_{i(1)} = 1|s_0) = P(Y_{i(2)} = 1|s_0) = 0.5$ and $P(Y_{i(1)} = 1, Y_{i(2)} = 0|s_0) = 0.25 = 0.5 \times 0.5$ by (A.3) in the Appendix. Here the sample variables $\{Y_{i(1)}, Y_{i(2)}\}$ under the posterior distribution given $s_0$ inherit the independence of the $Y$'s, even if the design were SRSWR.

If we selected $s_0 = \{1,1\}$, then $I_1^1(s_0)I_1^2(s_0) = 1$, $I_1^1(s_0)I_2^2(s_0) = I_2^1(s_0)I_1^2(s_0) = I_2^1(s_0)I_2^2(s_0) = 0$ and $P(Y_{i(1)} = 1, Y_{i(2)} = 0|s_0) = 0$. Here $\{Y_{i(1)}, Y_{i(2)}\}$ are stochastically dependent under $P(\cdot|s_0)$.

We next deal with a sequence of super-populations indexed by $\nu = 1, 2, \ldots$.

PROPOSITION 4.3 [Asymptotic normality under $P(\cdot|s_\nu)$]. *Let $Y_{\nu i}$, $i = 1, \ldots, N_\nu$, $\nu \geq 1$, be i.i.d. r.v.s on $(\Omega, \mathscr{F}, P)$ with zero mean and finite variance $\sigma^2 > 0$. Consider SRSWOR samples $s_\nu$ of size $n_\nu$ and $P_\nu = P(\cdot|s_\nu)$ as in (4.4). Let $Y_{\nu i(k)}(s_\nu, \omega)$ denote the array of r.v.s as in Definition 4.4. Then $(\sigma^2 n_\nu)^{-1/2}[\sum_{k=1}^{n_\nu} Y_{\nu i(k)}]$ converges in law to a standard normal random variable.*

The proof is in the Appendix.

PROPOSITION 4.4. *Let $B = \bigcup_{s \in A} \{s\} \times F_s \in C(S) \times \mathscr{F}$ with all $s$ distinct. We write*

$$(4.5) \qquad P_{d,m}(B|S \times \mathscr{F})(\omega) = \sum_{s \in A} p(s, \omega) I_{F_s}(\omega), \qquad \omega \in \Omega.$$

*Then the right hand-side of (4.5) is the conditional probability measure on $(S \times \Omega, C(S) \times \mathscr{F})$ given the $\sigma$-field $S \times \mathscr{F}$. The result is also valid if we replace everywhere $\mathscr{F}$ by $\mathscr{F}_N = \sigma(Y^N, X^N, Z^N)$ or by $\sigma(Z^N)$.*

An outline of the proof is given in the Appendix.

**5. Convergence in the product space and asymptotic independence.** In this section we establish results that enable us to determine the limiting distribution of a combination of sample estimators and super-population statistics. Let $\hat{\theta} \in \mathbb{R}^\ell$ be a sample estimator as in Remark 4.2. We define

$$F(t, \omega) = p(\{s \in S : \hat{\theta}(s, \omega) \leq t\}, \omega), \qquad t \in \mathbb{R}^\ell.$$

THEOREM 5.1. *We consider a sequence of product spaces and sample estimators as in Definition 4.3 and Remark 4.2, indexed by $\nu \geq 1$. Let $\lambda_\nu$, $\lambda \in \mathbb{R}^\ell$ be random vectors defined on $(\Omega, \mathscr{F}, P)$. We have:*

*(i) If $\lambda_\nu \to \lambda$ in the law of the model (P), then $\lambda_\nu \to \lambda$ in the law of the product space.*

*(ii) If $F_\nu(t, \omega) \to F(t, \omega)$ in probability $P$ for all points of continuity $t \in \mathbb{R}^\ell$ of $F(t, \omega)$, then $F(t, \omega)$ is a bounded random variable in the model space, and the product-space distribution of $\hat{\theta}_\nu$ converges to $F(t) = \int_\Omega F(t, \omega) \, dP(\omega)$. In particular, if $\hat{\theta}_\nu(\cdot, \omega)$ is design-consistent a.s. $\omega$, then it is consistent in the product space.*

*(iii) Assume that $\lambda_\nu \to \lambda$ in the law of the model and $F_\nu(t, \omega) \to F(t)$ in probability $P$ as $\nu \to \infty$ for all points of continuity $t \in \mathbb{R}_\ell$ of $F(t)$, where $F(t)$ is a nonstochastic distribution function. Then the joint distribution function of $(\hat{\theta}_\nu, \lambda_\nu)$ converges to the product of the two limiting distributions. The random variables $\hat{\theta}_\nu$ and $\lambda_\nu$ are said to be asymptotically independent.*

The proof is given in the Appendix. Note that when the limiting design-based distribution is normal with mean zero, we only require that the limiting variance be nonstochastic in the model. This last condition would follow if we imposed simple conditions in the super-population model, as we did in Proposition 3.1.

REMARK 5.1. The design-based distribution of the sample estimator $\hat{\theta}_N$ (viewed as a random variable in the product space) is a version of its conditional distribution in the product space given $S \times \mathscr{F}_N$. This follows if we take sets of the form $B(t) = \{(s,\omega) : \hat{\theta}(s,\omega) \leq t\}$, $t \in \mathbb{R}^\ell$, in Remark 4.1 and use (A.5) in the Appendix.

EXAMPLE 5.1 (The ratio estimator of the finite population mean). We assume a one-stage super-population model composed of $L$ disjoint strata of $N_h$ i.i.d. r.v.s $(Y_{hi}, Z_{hi})$, $i = 1, \ldots, N_h$, with mean $\mu_h = E_m(Y_{h1})$ and variance $\sigma_h^2 = V_m(Y_{h1})$, $h = 1, \ldots, L$. Let $\mu_N = \frac{1}{N} \sum_{h=1}^L N_h \mu_h$ be the parameter of interest, $N = N_1 + \cdots + N_L$ and $\Gamma_N = \frac{1}{N} \sum_{h=1}^L N_h \sigma_h^2$. The finite population mean is

$$\overline{Y}_N = \frac{1}{N} \sum_{h=1}^L \sum_{i=1}^{N_h} Y_{hi}.$$

Consider a stratified one-stage PPSWR design with the notation of Example 2.1. At each draw $k = 1, \ldots, n_h$, the units are selected in the sample $s_h$ with probabilities $p_{hi}$, which are functions of $Z_{hi}$, $i = 1, \ldots, N_h$, $h = 1, \ldots, L$. The ratio estimator of the finite population mean is

$$\bar{y}_R = (1/\widehat{N}) \sum_{h=1}^L \sum_{i \in s_h} y_{hi}/n_h p_{hi}, \qquad \widehat{N} = \sum_{h=1}^L \sum_{i \in s_h} 1/n_h p_{hi}.$$

Let $n = n_1 + \cdots + n_L$, $n_h \geq 1$ and $N \to \infty$, $n \to \infty$. We aim to obtain the asymptotic normality of $\sqrt{n}(\bar{y}_R - \mu_N)$ as $N \to \infty$. Here we construct a product space with the unconditional model probability measure $P$ rather than $P_z$ (defined after Definition 4.3). We decompose $\sqrt{n}(\bar{y}_R - \mu_N)$ into two terms, as in (1.1), and apply Theorem 5.1. The CLT for $\sqrt{N}(\overline{Y}_N - \mu_N)$ with limiting variance $\Gamma_m = \lim_N \frac{1}{N} \sum_{h=1}^L N_h \sigma_h^2$, $\Gamma_m < \infty$, follows if we assume Liapunov's condition (Theorem 27.3 in [2]),

$$\sum_{h=1}^L N_h E_m |Y_{h1} - \mu_h|^{2+\delta} = o(N^{1+\delta/2} \Gamma_m^{1+\delta/2}) \qquad \text{as } N \to \infty, \text{ for some } \delta > 0.$$

Let $\Gamma_d$, the limiting design variance of $\sqrt{n}(\bar{y}_R - \overline{Y}_N)$, that is, $\Gamma_d = \lim_n (1/N) \times (n/N) \sum_h (\sum_i e_{hi}^2(\omega)/n_h p_{hi} - e_h^2(\omega)/n_h)$, be positive definite, where $e_{hi}(\omega) = y_{hi}(\omega) - \overline{Y}_N(\omega)$ and $e_h(\omega) = \sum_i e_{hi}(\omega)$ are the residuals, $i = 1, \ldots, N_h$, $h =$ $1, \ldots, L$. Note that $C_2$ implies that $\widehat{N}$ is consistent. The CLT for $\sqrt{n}(\bar{y}_R - \overline{Y}_N)$ with asymptotic variance $\Gamma_d$ follows as in [31] by assuming conditions $(C_1)$ to $(C_3)$ in the Appendix applied to the residuals of a first stage sampling design, where $M = N$, and by Slutsky's theorem. Theorem 5.1 can then be applied if we assume that $\Gamma_d$ is nonstochastic.

**6. Sample estimators derived from an estimating equation (EE).** In this section we describe a methodology to derive the asymptotic normality of the root $\hat{\theta}_N \in \mathbb{R}^\ell$ of the sample estimating equation when centred about the super-population parameter $\theta_0 \in \mathbb{R}^\ell$. We combine existing asymptotic results in both the design and super-population probability as in Theorem 5.1.

Let $(\Omega, \mathscr{F}, P)$ and $(Y^{N\nu}, X^{N\nu}, Z^{N\nu})$ represent a super-population as in Definition 3.1 associated with a design space as in Definition 2.3. The first stage sample size is denoted by $n_\nu$. In what follows we omit the index $\nu$ and we set $N \to \infty$, $n \to \infty$ as $\nu \to \infty$. We first define a finite population estimating equation (EE) estimator and then an EE for the sample space.

DEFINITION 6.1. Let $g$ represent a continuously differentiable function defined on $\mathbb{R}^{p+k+\ell}$. We consider functions of the form

$$(6.1) \qquad G_N(\theta, \omega) = [1/\alpha(N)] \sum_{i=1}^N g(Y_i(\omega), X_i(\omega), \theta),$$

where $\omega \in \Omega$, $\theta \in \mathbb{R}^\ell$, $g \in \mathbb{R}^\ell$, $\alpha(N)/N = 0(1)$ as $N \to \infty$. A finite population EE is defined by

$$(6.2) \qquad G_N(\theta, \omega) = 0.$$

A finite population EE estimator is defined as a solution $\theta_N$ of (6.2), when such a solution exists and is unique. For $\omega \in \Omega$ fixed, $\theta_N$ is a finite population parameter.

DEFINITION 6.2. Let $\widehat{G}_N(\theta, \omega)$ be a design-consistent estimator of $G_N(\theta, \omega)$. A sample EE is defined by

$$(6.3) \qquad \widehat{G}_N(\theta, \omega) = 0.$$

A sample EE estimator $\hat{\theta}_N$ is defined as a solution of the sample EE in (6.3).

Yuan and Jennrich [30] (see also [3]) set general conditions for the existence, strong consistency and asymptotic normality of EE estimators which require independent but not necessarily i.i.d. random vectors $g(Y_i, X_i, \theta)$, $i = 1, \ldots, N$. We can apply their results to clustered data models with cluster totals $g_i(\theta) = \sum_{j=1}^{M_i} g(Y_{ij}, X_{ij}, \theta)$, which are stochastically independent. The cluster sizes $M_i$, $i = 1, \ldots, N$, stay bounded as the number $N$ of clusters goes to infinity. Theorem 6.1 shows that the sample EE estimator (around the model parameter) is asymptotically normal in the law of the product space. Conditions 1–3 were given by Yuan and Jennrich [30] for the existence and consistency of $\theta_N$ and the asymptotic normality of $\sqrt{N}(\theta_N - \theta_0)$. Conditions 1, 4 and 5 below imply the existence and design-consistency of $\hat{\theta}_N$ and the design-asymptotic normality of $\sqrt{n}(\hat{\theta}_N - \theta_N)$.

THEOREM 6.1. *Consider a sequence of super-populations composed of $N$ independent random vectors associated to design spaces as in Definition 3.1. Let $f = \lim_n n/N \geq 0$ as $n \to \infty$. Note that we do not require that $f = 0$. We assume the following conditions.*

1. $G_N(\theta_0) \to 0$ *with probability one.*
2. *There is a compact neighborhood $B(\theta_0)$ of $\theta_0$ on which, with probability one, all $G_N(\theta)$ are continuously differentiable and the Jacobians $\partial G_N(\theta)/\partial \theta$ converge uniformly in $\theta$ to a nonstochastic limit $J(\theta)$ which is nonsingular at $\theta_0$.*
3. $\sqrt{N} G_N(\theta_0) \Rightarrow N(0, \Gamma_m)$ *in the law of the super-population.*
4. *There is a compact neighborhood $B(\theta_0)$ of $\theta_0$ on which $\partial \widehat{G}_N(\theta)/\partial \theta$ converge uniformly in the design probability to a nonstochastic (in design) limit which coincides with $J(\theta)$ at $\theta_0$ for almost every $\omega \in \Omega$. Note that if $\widehat{G}_N(\theta) = \sum_{i \in s} c_i g_i(\theta)$, $c_i$ independent of $\theta$ and $\{G_N(\theta), N \geq 1\}$ are continuously differentiable, then $\{\widehat{G}_N(\theta), N \geq 1\}$ are also continuously differentiable.*
5. $\sqrt{n}\widehat{G}_N(\theta_N) \Rightarrow N(0, \Gamma_d)$ *in the law of the design as $n \to \infty$ for almost every $\omega \in \Omega$, where the variance–covariance $\Gamma_d$ is nonstochastic in the super-population. Let $J = J(\theta_0)$ and $\Gamma = J^{-1}[\Gamma_d + f\Gamma_m]J^{-1}$. Then we have, in the law of the product space,*

(6.4) $$\sqrt{n}(\hat{\theta}_N - \theta_0) \Rightarrow N(0, \Gamma).$$

Estimation of $\Gamma$ from the sample data depends on the particular design under consideration for the estimation of $\Gamma_d$, and on both the model assumed for the variance–covariance structure of the super-population and the sampling design for the estimation of $\Gamma_m$ (see Example 6.1). The Jacobian matrix $J = J(\theta_0)$ can be estimated consistently by $(\partial \widehat{G}_N/\partial \theta)(\hat{\theta}_N)$: this follows from Assumptions 2 and 4 and the consistency of $\hat{\theta}_N$ [from (6.4)]. The proof of Theorem 6.1 is in the Appendix.

REMARK 6.1. Korn and Graubard [14] propose direct estimators of the variance–covariance of the sample mean under different super-population models and sampling designs. See also Rubin–Bleuer [21]. In Example 6.1 we assume a two-stage super-population model and design to estimate $\Gamma$.

In Example 6.1, (i) we establish a model for the super-population variance-covariance structure so that we can estimate the model variance matrix $\Gamma_m$ from the sample data and (ii) we examine the design conditions for the asymptotic normality of $\sqrt{n}\widehat{G}_N(\theta_N)$ to hold in the finite population.

EXAMPLE 6.1 (General EE sample estimator under a stratified two-stage super-population model and design). Assume the stratified two stage super-population model of Example 3.1 with the addition of the auxiliary information given by $X_{hij}$ ($h, i$ and $j$ as in Example 3.1) and the two-stage design of Example 2.1. In this example, we construct separate consistent estimators within each stratum and so we need to assume that $n_h \to \infty$ and $N_h \to \infty$, $h = 1, \ldots, L$.

Let the finite population EE, given by

$$G_N(\theta) = \frac{1}{M} \sum_{h=1}^{L} \sum_{i=1}^{N_h} g_{hi}(\theta),$$

with $g_{hi}(\theta) = \sum_{j=1}^{M_{hi}} g(Y_{hij}, X_{hij}, \theta)$, satisfy the first three conditions of Theorem 6.1. Now let $\widehat{G}_N(\theta)$ be the sample estimator of $G_N(\theta)$, where $\widehat{G}_N(\theta)$ replaces $\hat{\theta}_N$ in Example 2.1. Assume $M/N \to \mathfrak{m} < \infty$ as $N \to \infty$. Also assume $n_h/N_h = c_h$ constant as $N_h \to \infty$, for all $h = 1, \ldots, L$.

(i) Assume that the second stage observations $g_{hij}(\theta)$ are i.i.d. r.v.s with means $\mu_{hi}$ and variances $\sigma_{hi}^2$, $j = 1, 2, \ldots, M_{hi}$. Furthermore, $(\mu_{hi}, \sigma_{hi}^2)$ are i.i.d. r.v.s, where the $\mu_{hi}$s have model variances $V_m(\mu_{hi}) = \gamma_h$, and the $\sigma_{hi}^2$s have model expectations $E_m(\sigma_{hi}^2) = \sigma_h^2$, $i = 1, 2, \ldots, N_h$, $h = 1, \ldots, L$. Thus,

$$V_m(\sqrt{N}G_N(\theta)) = (N/M) \sum_{h=1}^{L} \left\{ W_h \sigma_h^2 + \gamma_h \left( \sum_{i=1}^{N_h} M_{hi}^2/M_h \right) \right\}$$

(6.5)

$$\text{with } W_h = \sum_{i=1}^{N_h} M_{hi}/M, \ h = 1, \ldots, L.$$

To obtain a (model) consistent estimator of $V_m(\sqrt{N}G_N(\theta))$, it is enough to get consistent or asymptotically unbiased estimators of $\sigma_h^2$ and $\gamma_h$, $h = 1, \ldots, L$. These can be written as quadratic functions of the finite population values $g_{hij} = g_{hij}(\theta)$ and, thus, they are finite population parameters:

$$\tilde{\sigma}_h^2 = \frac{1}{N_h} \sum_{i=1}^{N_h} \left\{ \sum_{j=1}^{M_{hi}} g_{hij}^2 - g_{hi}^2/M_{hi} \right\} \Big/ (M_{hi} - 1), \qquad h = 1, \ldots, L,$$

and

$$\tilde{\gamma}_h = \frac{1}{N_h - 1} \left\{ \sum_{i=1}^{N_h} \left\{ \frac{g_{hi}}{M_{hi}} \right\}^2 - \frac{1}{N_h} \left( \sum_{i=1}^{N_h} g_{hi}/M_{hi} \right)^2 \right\}, \qquad h = 1, \ldots, L,$$

where $\tilde{\sigma}_h^2$ and $\tilde{\gamma}_h$, $h=1,\ldots,L$, are model unbiased, as well as model consistent. Any pair of design-consistent or asymptotically design-unbiased estimators $\hat{\sigma}_h^2$ and $\hat{\gamma}_h$ of the respective finite population parameters $\tilde{\sigma}_h^2$ and $\tilde{\gamma}_h$ can replace $\sigma_h^2$ and $\gamma_h$, $h=1,\ldots,L$, in (6.5) to yield an asymptotically unbiased estimator of $\Gamma_m = \lim_N V_m(\sqrt{N} G_N(\theta))$ in the product space.

(ii) To obtain the asymptotic normality of $\sqrt{n}\widehat{G}_N(\theta_N)$, we express $\widehat{G}_N(\theta_N)$ as the sum of $n = n_1 + \cdots + n_L$ independent, zero mean random vectors $Z_{hk}(\theta_N)$,

$$\widehat{G}_N(\theta_N) = \hat{G}_N(\theta_N) - G_N(\theta_N), \qquad \widehat{G}_N(\theta_N) = \sum_{h=1}^L W_h \frac{1}{n_h} \sum_{k=1}^{n_h} Z_{hk}(\theta_N)$$

with $W_h = M_h/M$, $Z_{hk}(\theta_N) = \sum_{i=1}^{N_h} \hat{g}_{hi}(\theta_N) I_{hi}^k / M_{hi} - \sum_{i=1}^{N_h} g_{hi}(\theta_N)/M_h$, where $\hat{g}_{hi}(\theta_N)$ denotes the second stage unbiased sample estimator of $g_{hi}(\theta_N)$.

The design-based CLT for $\sqrt{n}\widehat{G}_N(\theta_N)$ with positive definite $\Gamma_d = \lim_n n \times \sum_{h=1}^L W_h^2 V_d(Z_{k1}(\theta_N))/n_h$ follows from conditions ($C_1$)–($C_3$) in the Appendix, with $\hat{\theta}_h^k$ replaced by $Z_{hk}(\theta_N)$.

As in Proposition 3.1, one can give conditions in the super-population so that a Liapunov-type condition holds in the design space and the asymptotic design variance $\Gamma_d$ exists and is nonstochastic (condition 5 of Theorem 6.1). The super-population conditions required for the latter are more complex than those stated in Proposition 3.1, but they can be specified in the same way. We do not spell them out here.

## APPENDIX

**Yung and Rao [31] designs conditions for the asymptotic normality of the sample mean.**

($C_1$) $n^{1+\delta} \sum_{h=1}^L \sum_{k=1}^{n_h} E_d |W_h \hat{\theta}_h^k / n_h|^{2+\delta} = O(1)$ as $n \to \infty$, $\hat{\theta}_h^k$ as in Example 2.1.

($C_2$) $(n/M) \max_{h,i,j} m_{hi} w_{hij} = O(1)$ as $n \to \infty$, where $m_{hi}$ are the second stage sample sizes and $w_{hij}$ are the sampling weights.

($C_3$) $\Gamma_d^N(\omega) = n V_d(\hat{\theta}_N) \to \Gamma_d$ positive definite as $n \to \infty$, $\hat{\theta}_N$ as in Example 2.1.

PROOF OF PROPOSITION 3.1. Since if $E_m|X|$ is finite, then $|X(\omega)|$ is finite a.s. $\omega$, condition $M_1$ implies $(1/N) \sum_{h=1}^L \sum_{i=1}^{N_h} |Y_{hi}|^{2+\delta} = O(1)$ a.s. $\omega$. ($C_1'$) follows from the boundedness of two terms once we take $N=2$ and $p=2+\delta$ in the inequality $E_d|(1/N)\sum_{k=1}^N X_k|^p \le (1/N)\sum_{k=1}^N E_d|X_k|^p$ (see (7), page 95 of [5]). Since $|\overline{Y}_h|^{2+\delta} = |E_d[\hat{\theta}_h^k]|^{2+\delta} \le E_d|\hat{\theta}_h^k|^{2+\delta}$, we only need to show that, for all $k=1,\ldots,n_h$, $h=1,\ldots,L$, $\sum_{h=1}^L W_h E_d|\hat{\theta}_h^k|^{2+\delta} =$

$O(1)$ a.s. $\omega$ as $N \to \infty$. At the $k$-draw we select one cluster, so $E_d|\hat{\theta}_h^k|^{2+\delta} = \sum_{i=1}^{N_h} |Y_{hi}(\omega)/M_{hi}|^{2+\delta} p_{hi} = \sum_{i=1}^{N_h} |Y_{hi}(\omega)|^{2+\delta} M_{hi}^{-1-\delta} M_h^{-1}$. Since $N \leq M$, $\sum_{h=1}^{L} W_h E_d|\hat{\theta}_h^k|^{2+\delta} \leq (1/N) \sum_{h=1}^{L} \sum_{i=1}^{N_h} |Y_{hi}|^{2+\delta}$, which is $O(1)$ when $(M_1)$ holds. $\square$

PROOF OF REMARK 4.3. Under SRS we have $P_{d,m}(Y_{i(k)}(s,\omega) \leq a) = (1/N) \sum_{i=1}^{N} P(Y_i(\omega) \leq a)$, $k = 1, \ldots, n$. For $n \geq 2$, $k \neq \ell$, $k, \ell = 1, \ldots, n$, under SRSWOR we have

$$P_{d,m}(Y_{i(k)}(s,\omega) \leq a, \ Y_{i(\ell)}(s,\omega) \leq b)$$
$$= 1/[N(N-1)] \sum_{i} \sum_{i \neq j} P(Y_i(\omega) \leq a) P(Y_j(\omega) \leq b), \tag{A.1}$$

and under SRSWR we have

$$P_{d,m}(Y_{i(k)}(s,\omega) \leq a, Y_{i(\ell)}(s,\omega) \leq b)$$
$$= (1/N^2) \bigg\{ \sum_{i} \sum_{i \neq j} P(Y_i(\omega) \leq a) P(Y_j(\omega) \leq b) + \sum_{i} P(Y_i(\omega) \leq \min(a,b)) \bigg\}. \tag{A.2}$$

Under SRSWOR let $P(Y_i \leq a) = p(a)$ for all $i = 1, \ldots, N$. Then $P_{d,m}(Y_{i(k)}(s,\omega) \leq a) = p(a)$, $k = 1, \ldots, n$, and the right-hand side in (A.1) is $p(a)p(b)$, which proves pairwise independence in the product space. Overall independence is proved similarly. If the $Y$'s are not identically distributed, we show dependence via a counterexample. Let $P(Y_1 \leq a) = p_1$ and $P(Y_i \leq a) = p_2$ for $i = 2, 3, \ldots, N$. If we take $N = 2$, we have $P(Y_{i(k)} \leq a) = [p_1 + p_2]/2$, $k = 1, 2$, and $P(Y_{i(1)} \leq a, Y_{i(2)} \leq a) = p_1 p_2$. Independence holds only when $p_1 = p_2$. Under SRSWR dependence in the product space follows from (A.2). For discrete $Y$'s,

$$P_{d,m}(Y_{i(k)}(s,\omega) = a, \ Y_{i(\ell)}(s,\omega) = b)$$
$$= (1/N^2) \bigg\{ \sum_{i} \sum_{i \neq j} P(Y_i(\omega) = a) P(Y_j(\omega) = b) \bigg\}, \qquad a \neq b. \tag{A.2'}$$

$\square$

PROOF OF PROPOSITION 4.2. For $s_0 \in S$, $k \neq \ell$, $k, \ell = 1, \ldots, n$, we have, by Proposition 4.1, part (ii),

$$P(Y_{i(k)}(s,\omega) \leq a|s_0) = \sum_{i=1}^{N} P(Y_i \leq a) I_i^k(s_0)$$

and

$$P(Y_{i(k)}(s,\omega) \leq a, \ Y_{i(\ell)}(s,\omega) \leq b|s_0)$$

(A.3)
$$= \sum_{i=1}^{N} \sum_{j=1}^{N} P(Y_i \leq a, \ Y_j \leq b) I_i^k(s_0) I_j^\ell(s_0).$$

In the WOR case we have $I_i^k(s_0)I_i^\ell(s_0) = 0$ for every $k \neq \ell$, $i = 1, \ldots, N$, and these terms disappear in the double sum, yielding independence. For samples $s_0 \in S$ for which $I_i^k(s_0)I_i^\ell(s_0) = 1$ for some $i$, the double sum above contains nonzero terms where $i = j$. The terms corresponding to the repeated labels in the product of the two distributions are different from their counterpart terms in the joint distribution:

$$\sum_{i=1}^{N} P(Y_i \leq a) P(Y_i \leq b) I_i^k(s_0) I_i^\ell(s_0)$$

$$\neq \sum_{i=1}^{N} P(Y_i \leq \min(a,b)) I_i^k(s_0) I_i^\ell(s_0) \qquad \text{for continuous } Y\text{'s}$$

and

$$\sum_{i=1}^{N} P(Y_i = a) P(Y_i = b) I_k^i(s_0) I_i^\ell(s_0)$$

$$\neq \sum_{i=1}^{N} P(Y_i = a, Y_i = b) I_i^k(s_0) I_i^\ell(s_0) \qquad \text{for discrete } Y\text{'s}. \qquad \square$$

PROOF OF PROPOSITION 4.3. Under $P(\cdot|s_\nu)$ the $Y_{\nu i(k)}$, $i = 1, \ldots, N_\nu$, $\nu \geq 1$, are i.i.d. r.v.s with mean zero and constant variance. That they are identically distributed like the original $Y$'s follows by (A.3), and independence follows from Proposition 4.2 for SRSWOR. As in Theorem 27.2 of [2], (27.9) holds, which implies the Lindeberg condition and proves the result. $\square$

PROOF OF PROPOSITION 4.4. We prove that, for each $\omega \in \Omega$, the right-hand side of (4.5) is a probability measure on the product space, and that, for each measurable set $B$ in the product space, it is a version of the conditional probability of $B$ given $S \times \mathscr{F}$ ([5], page 223). The first statement follows from the additivity of $p$ and the $\sigma$-additivity of the indicator functions. To prove the second, we note first that $p(s, \cdot)$ is $\mathscr{F}$-measurable. Then it suffices to show that, on the elementary rectangles $B = s_0 \times F_0$, we have

$$\int_{S \times F} p(s_0, \omega) I_{F_0}(\omega) \, dP_{d,m} = P_{d,m}(B \cap (S \times F)), \qquad F \in \mathscr{F}.$$

The left-hand side above equals $\sum_{s\in S} \int_{F\cap F_0} p(s_0,\omega)p(s,\omega)\,dP$ by (4.2).

By Definition 2.3(ii), the sum above equals

$$\int_{F\cap F_0} p(s_0,\omega)\,dP = P_{d,m}(\{s_0\}\times F\cap F_0). \qquad \square$$

PROOF OF THEOREM 5.1. Let $t\in\mathbb{R}^\ell$. We first omit indexing the populations. Let

$$B(t)=\{(s,\omega):\hat{\theta}(s,\omega)\le t\}=\bigcup_{s\in S}\{s\}\times F_s$$

with

$$F_s=\{\omega\in\Omega:\hat{\theta}(s,\omega)\le t\}\in\mathscr{F}.$$

By Remark 4.1, $P_{d,m}(B(t))=\sum_{s\in S}\int_\Omega p(s,\omega)I_{F_s}(\omega)\,dP$. Note that

$$(A.4) \qquad F(t,\omega)=p(\{s:\hat{\theta}(s,\omega)\le t\},\omega)=\sum_{s\in S} p(s,\omega)I_{A_\omega}(s),$$

where $A_\omega=\{s\in S:\hat{\theta}(s,\omega)\le t\}\in C(S)$. For each $(s,\omega)$ the indicator function of $F_s$ coincides with the indicator function of $A_\omega$,

$$(A.5) \qquad I_{F_s}(\omega)=I_{A_\omega}(s).$$

Using (A.5) in (A.4) and the formula for $P_{d,m}(B(t))$, we have

$$(A.6) \quad \begin{aligned} P_{d,m}(B(t))&=\int_\Omega F(t,\omega)\,dP \quad\text{and}\\ P_{d,m}(B(t)\cap E)&=\int_{\Omega\cap E} F(t,\omega)\,dP, \qquad E\in\mathscr{F}.\end{aligned}$$

(i) This follows from $P(\lambda_\nu\le u)=P_{d,m}(S\times\{\lambda_\nu\le u\})$ and $\sum_{s\in S} p(s,\omega)=1$ for all $\omega\in\Omega$.

(ii) $F_\nu(t,\omega)$ converges in probability to $F(t,\omega)$ at points of continuity $t$. Since $0\le F_\nu(t,\omega)\le 1$, the bounded convergence theorem (Theorem 16.5 of [1], page 180) implies

$$(A.7) \qquad \int_\Omega F_\nu(t,\omega)\,dP - F(t) \to 0 \qquad \text{as } \nu\to\infty.$$

(A.6) and (A.7) yield (ii),

$$P^\nu_{d,m}(B_\nu(t))=\int_\Omega F_\nu(t,\omega)\,dP \to F(t) \qquad \text{as } \nu\to\infty.$$

(iii) Consider the indicator function $I_\nu(u,\omega) = I_{\{\omega : \lambda_\nu(\omega) \le u\}}(\omega)$. Using (A.6) with $E_\nu(u) = S_\nu \times \{\omega : \lambda_\nu(\omega) \le u\}$, the joint distribution function of $(\hat{\theta}_\nu, \lambda_\nu)$ can be expressed as

$$P_{d,m}^\nu \{(s,\omega) : (\hat{\theta}_\nu, \lambda_\nu) \le (t,u)\} = P_{d,m}(B_\nu(t) \times E_\nu(u)) = \int_\Omega I_\nu(u,\omega) F_\nu(t,\omega)\, dP,$$

and if $H$ denotes the distribution function of $\lambda$, we have, for points of continuity $t$, $u$ of $F$ and $H$,

$$\int_\Omega I_\nu(u,\omega) \cdot F_\nu(t,\omega)\, dP - F(t)H(u)$$
$$= \int_\Omega I_\nu(u,\omega)(F_\nu(t,\omega) - F(t))\, dP + F(t) \int_\Omega (I_\nu(u,\omega) - H(u))\, dP.$$

All functions are bounded by one, so the first term of the right-hand side converges to zero by the bounded convergence theorem since, by hypothesis, $F_\nu(t,\cdot) - F(t)$ converges to zero in probability $P$ at points of continuity $t$ of $F$. The second term also converges to zero by hypothesis. □

PROOF OF THEOREM 6.1. For simplicity we assume that $f = n/N$ for all $n$:

(A.8) $\qquad \sqrt{n}(\hat{\theta}_N - \theta_0) = \sqrt{n}(\hat{\theta}_N - \theta_N) + \sqrt{f}\sqrt{N}(\theta_N - \theta_0).$

Assumptions 1–3 imply the asymptotic normality of the second term on the right-hand side of (A.8), in the law of the model (see [30]). This and Theorem 5.1(i) imply convergence in the law of the product space. Next we observe that $\hat{\theta}_N$ exists and $\hat{\theta}_N - \theta_N \to 0$ in design probability as $n \to \infty$ for almost every $\omega \in \Omega$. Indeed, the $\widehat{G}_N(\theta)$ are continuously differentiable and design consistency implies that $\widehat{G}_N(\theta)$ converges to $G(\theta)$ (the limit of $G_N(\theta)$ in [30]) in design probability. Hence, we can apply to $\widehat{G}_N(\theta)$ the techniques of Theorems 1 and 2 of [30], and, thus, Assumptions 1 and 4 imply that $\hat{\theta}_N \to \theta_0$ in the design probability. Since the above mentioned theorems also imply $\theta_N \to \theta_0$ a.s. in the model probability $P$, we have $\hat{\theta}_N - \theta_N \to 0$ in the design probability a.s. $\omega$. Conditions 4 and 5 imply asymptotic normality of the first term in the right-hand side of (A.8). This, in turn, implies convergence in the product space, by Theorem 5.1(ii). The two terms in (A.8) are not stochastically independent in general. Theorem 5.1(iii) and Assumption 5 imply the "asymptotic independence" of the terms and the asymptotic normality of the sum. □

**Acknowledgments.** The authors would like to thank J. N. K. Rao for his advice and guidance, and for his continuing support and encouragement. We are indebted to the referee for his countless helpful suggestions and for a very thorough revision which resulted in a much more accurate and readable paper. In addition, we would like to thank the referee for providing a numerical illustration which serves as a basis for our Examples 4.1 and 4.3. And last but not least, we would like to thank Jack Gambino for his very useful comments.

BUSINESS SURVEY METHODS DIVISION
STATISTICS CANADA
11TH FLOOR, RH COATS BUILDING
TUNNEYS PASTURE
OTTAWA, ONTARIO
CANADA K1A 0T6
E-MAIL: rubisus@statcan.ca